\documentclass{amsart}

\usepackage{amsmath}
\usepackage{graphicx}
\usepackage{url}
\usepackage{mathrsfs}

\newtheorem*{theorem}{Theorem}
\newtheorem*{conjecture}{Conjecture}
\newtheorem{algorithm}{Algorithm}
\newtheorem*{definition}{Definition}

\theoremstyle{definition}
\newtheorem{move}{Move}

\theoremstyle{remark}
\newtheorem*{remark}{Remark}

\begin{document}

\title[Computer analysis of Sprouts with nimbers]{Computer analysis of Sprouts with nimbers}

\author{Julien Lemoine - Simon Viennot}

%\date{\today}

\begin{abstract}
Sprouts is a two-player topological game, invented in 1967 in the University of Cambridge by John Conway and Michael Paterson. The game starts with $p$ spots, and ends in at most $3p-1$ moves. The first player who cannot play loses.

The complexity of the $p$-spot game is very high, so that the best hand-checked proof only shows who the winner is for the $7$-spot game, and the best previous computer analysis reached $p=11$.

We have written a computer program, using mainly two new ideas. The nimber (also known as Sprague-Grundy number) allows us to compute separately independent subgames; and when the exploration of a part of the game tree seems to be too difficult, we can manually force the program to search elsewhere. Thanks to these improvements, we reached up to $p=32$. The outcome of the 33-spot game is still unknown, but the biggest computed value is the 47-spot game ! All the computed values support the \textit{Sprouts conjecture}: the first player has a winning strategy if and only if $p$ is $3$, $4$ or $5$ modulo $6$.

We have also used a check algorithm to reduce the number of positions needed to prove which player is the winner. It is now possible to hand-check all the games until $p=11$ in a reasonable amount of time.
\end{abstract}

\maketitle

\section{Introduction}

Sprouts is a two-player pencil-and-paper game invented in 1967 in the University of Cambridge by John Conway and Michael Paterson\cite{mg67}. The game starts with $p$ spots and players alternately connect the spots by drawing curves between them, adding a new spot on each curve drawn. A new curve cannot cross or touch any existing one, leading necessarily to a \emph{planar} graph. The first player who cannot play loses\footnote{The \emph{mis\`ere} version of the game corresponds to the opposite rule, i.e. a player unable to move wins. The analysis of the mis\`ere version is more difficult, and is the object of another article\cite{lv09}.}.\par
The last rule states that a spot cannot be connected to more than 3 curves, and it induces the end of the game in a finite number of moves. If we consider that a spot has 3 \emph{lives} at the beginning of a game, there is a total number of $3p$ lives. Each move then consumes 2 lives and creates a spot with 1 life, globally decreasing by one the total number of lives. It follows that the game ends in at most $3p-1$ moves.
\begin{figure}[ht]
\centering
\includegraphics[scale=0.5]{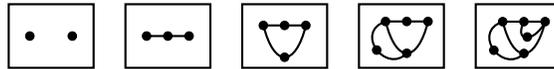}
\caption{A sample game of 2-spot Sprouts (the second player wins).}
\end{figure}
\par
Despite the small number of moves of a given game, it is difficult to determine whether the winning strategy is for the first or the second player. The best published and complete hand-checked proof is due to Focardi and Luccio, and shows who the winner is for the $7$-spot game\cite{fl04}.
\par
The first published computer analysis of Sprouts is due to Applegate, Jacobson and Sleator \cite{ajs91}: they computed in 1991 which player wins up to the 11-spot game. They noted a pattern in their results and proposed the \textit{Sprouts conjecture}: the first player has a winning strategy in the $p$-spot game if and only if $p$ is $3$, $4$ or $5$ modulo $6$. We have written a new program, using an improved version of their original representation of Sprouts positions, and obtained significant new results with the help of two main ideas: the theoretical concept of \emph{nimber} and the manual exploration of the game tree.
\par
Our program enabled us to compute which player wins the $p$-spot game up to 32 spots, and also some values up to 47 spots. All the computed values support the Sprouts conjecture. We then performed a second computation to validate the results, in which our program tries to minimize the number of positions needed to prove the value of a given $p$-spot game. It is then possible to generate graphs corresponding to these minimized databases, which provides a hand-checkable winning strategy for the $p$-spot game.

\section{Game tree}

We will call \textit{position} a graph embedded in the plane, obtained with the rules of Sprouts.\par

From any given Sprouts position, there is necessarily a winning strategy, either for the player to move, or for the next player. If the winning strategy is for the player to move, we will say that the \emph{outcome} of the position is a \emph{win} and in the other case, that it is a \emph{loss}. The main goal of our program is then to compute the outcome of Sprouts positions, and particularly the starting positions with $p$ spots.

\begin{figure}[ht]
\centering
\includegraphics[scale=0.5]{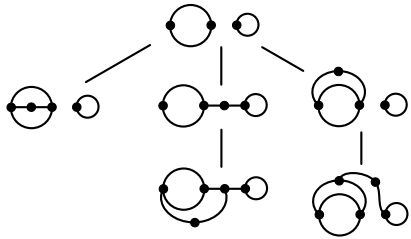}
\caption{Game tree of a given Sprouts position.}
\label{trees}
\end{figure}

We call \textit{game tree of a position} the tree where the \emph{root} is the given position, the nodes are the positions that can be reached by playing moves from the root, and where an edge links two positions if the \emph{child} is obtained from the \emph{parent} in only one move. Figure \ref{trees} is an example of such a tree.

As we will see in section \ref{mainalgo}, the outcome of a position is computed recursively by developing its game tree.

\section{Positions representation}
\label{posrepr}

As a pencil-and-paper game, Sprouts is not suitable for programming. We describe here a representation of the game with strings that could be used with a computer to develop game trees, and therefore deduce winning strategies.

\subsection{Regions and boundaries}

In a given position, a \textit{region} is a connected component of the plane and inside a given region, a \textit{boundary}, is a connected component of the curves drawn by the players. An isolated spot is considered as a boundary in itself.

\begin{figure}[ht]
\centering
\includegraphics[scale=0.5]{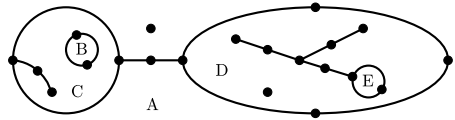}
\caption{A Sprouts position obtained after 10 moves in a 11-spot game.}
\label{game1}
\end{figure}

For example, figure \ref{game1} contains 5 regions and 9 boundaries: there are 3 boundaries in region D, 2 boundaries in regions A and C and 1 boundary in regions B and E.

\subsection{String representation}

In order to perform computations with a computer, we need a way to represent the positions of Sprouts by strings of characters. We follow the representation described in \cite{ajs91} in 1991\footnote{The equivalence between graphic and string representations is easy to see, but hard to demonstrate. This is the goal of a full 400 pages report\cite{dhkr90}.}. Basically, each spot (\emph{vertex} in the graph theory) will be denoted by a letter, and the graph is described as follows:
\begin{itemize}
 \item The complete graph is represented by the list of strings of the regions that it contains and is terminated by an \textit{end-of-position} character: ``\verb/!/''
 \item A region is represented by the list of strings of the boundaries that it contains and is terminated by an \textit{end-of-region} character: ``\verb/}/''
 \item A boundary is represented by the list of its vertices and is terminated by an \textit{end-of-boundary} character: ``\verb/./''
\end{itemize}

For a given boundary, we write the vertices in the order they appear while following one side of the boundary. Inside a given region, the same orientation must be respected.

In our example, we turn around all boundaries clockwise, except for the (unique) boundary that surrounds a region, around which we turn counterclockwise. A string representing the above position is then:

\begin{center}
 \verb/AL.}AL.BNMCMN.}D.COFPGQFOCM.}E.HRISJSIUKTKUIR.FQGP.}KT.}!/
\end{center}

\medskip
\begin{center}
\includegraphics[scale=0.5]{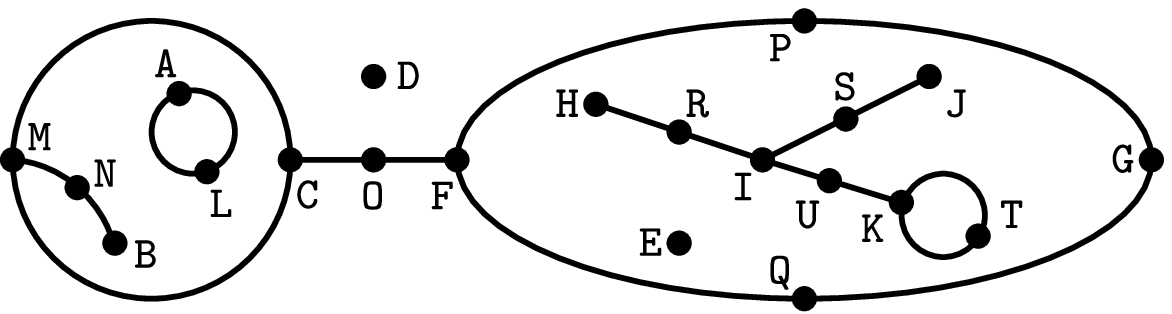}
\end{center}
\medskip

The boundary \verb/AL./ appears twice, because it is inside the same region as the 4 vertices \verb/B;N;M;C/, and it also surrounds a region itself. Notice how the string of the boundary containing \verb/B;N;M;C/ is obtained: starting at \verb/B/, we follow the bottom of the boundary, meeting first \verb/N/, then \verb/M/, and \verb/C/. Then, we continue to follow the top of the boundary, meeting \verb/M/ and \verb/N/ once again (but on the opposite side), stopping when we are back at the starting point.\par

\begin{remark}
The strings of the starting positions are \verb/A.}!/ (for the $1$-spot game), \verb/A.B.}!/ ($2$-spot game), \verb/A.B.C.}!/ ($3$-spot game)...\par
\end{remark}

\subsection{Equivalent positions and strings}

\begin{figure}[ht]
\centering
\includegraphics[scale=.5]{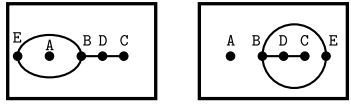}
\caption{Two equivalent positions.}
\label{2eqpos}
\end{figure}

A single position can correspond to several strings. For example, the first position in figure \ref{2eqpos}
 can lead to \verb/BDCDBE.}BE.A.}!/ or to \verb/A.EB.}DCDBEB.}!/. In the same way, a single string can represent several positions: \verb/BDCDBE.}BE.A.}!/ can represent the two positions of figure \ref{2eqpos}.

More precisely, a single string represents an infinite number of positions, as we can draw lines of multiple shapes as long as the topology is respected. Conversely, a single position corresponds to only a finite number of strings if we only allow the strings to use the first letters of the alphabet.

\subsection{Moves}

The main interest of the string representation described above is that moves can be easily defined and computed from it. Moves always take place in a single region but they are of two different types, depending on whether we link one boundary to itself or two different boundaries. The program we describe below will thus produce finite game trees.

\begin{move}
A \textit{two-boundaries} move consists in connecting two spots of different boundaries.

Let $\tt x \rm _1...\tt x \it _m$ and $\tt y \rm _1...\tt y \it _n$ be two different boundaries in the same region, with $m\geq 2$ and $n\geq 2$.
We suppose that $\tt x \it _i$ and $\tt y \it _j$ are vertices that occur two times or less in the whole string, with $1\leq i\leq m$ and $1\leq j\leq n$.\par

Then the two-boundaries move consists in merging these two boundaries in\\ $\tt x\rm_1...\tt x\it _i\tt zy\it _j...\tt y \it_n \tt y\rm _1...\tt y\it _j\tt zx\it _i...\tt x\it _m$ where $\tt z$ is the new created vertex.\par

The same definition holds if $m=1$ or if $n=1$, but in these cases, $\tt x\it _i...\tt x\it _m$ and $\tt y\it _j...\tt y\it _n$ are empty boundaries.
\end{move}

\begin{move}
A \textit{one-boundary} move consists in connecting two spots of the same boundary.

Let $\tt x\rm _1...\tt x\it _n$ be a boundary, with $n\geq 2$. We suppose that $\tt x\it _i$ and $\tt x\it _j$ are vertices that occur two times or less in the whole string, with $1\leq i\leq m$, $1\leq j\leq n$ and $i\neq j$, or that $i=j$ and $\tt x\it _i$ occurs only once in the whole string.\par

Then the one-boundary move consists in separating the other boundaries of the same region in 2 sets $\tt B\rm _1$ and $\tt B\rm _2$, and the original region is divided into two new regions: $\tt x\rm _1...\tt x\it _i\tt z\tt x\it _j...\tt x\it _n.\tt B\rm _1\}$ and $\tt x\it _i...\tt x\it _j\tt z.B\rm _2\}$.\par

The same definition holds if $n=1$, but in this case, $\tt x\it _j...\tt x\it _n$ is an empty boundary.
\end{move}

\begin{figure}[ht]
\centering
\includegraphics[scale=.5]{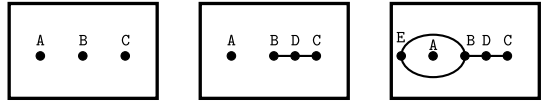}
\caption{Two moves in a 3-spot game.}
\label{3spotgame}
\end{figure}

Let us give an example of the two types of moves with a $3$-spot game. The initial string is \verb/A.B.C.}!/. First, we make a two-boundaries move, which leads to \verb/A.BDCD.}!/, and then we make a one-boundary move, which separates the position into two regions: \verb/BDCDBE.}BE.A.}!/. Figure \ref{3spotgame} shows three corresponding positions.

 These definitions are sufficient to create a first version of a program for computing Sprouts. This program could determine the outcome of the $p$-spot game for a few values of $p$, but there would be too many equivalent strings, and the memory would saturate quite quickly.

\section{Strings simplification}

We will now explain several methods to simplify the strings, with two main goals. Firstly, if we go back to figure \ref{trees}, the middle and the right children of the root can be represented by the same string (we will know at the end of this section that it is \verb/22.}]!/). These kind of simplifications will merge some equivalent branches of the game tree and thus decrease the complexity of the computation.

Secondly, the simplified strings will be more suited to perform the \emph{canonization} step described in section \ref{canosec}.

\subsection{Deletion of the dead parts}

First of all we delete the dead vertices (those which occur 3 times in the string). Then, we delete the empty boundaries (boundaries whose vertices are all dead) and finally, we delete the dead regions (with 0 or 1 life).

Using our previous example, we would therefore delete 5 dead vertices. Then, no boundary would be empty, but the region \verb/KT.}!/ would be dead. The new string and position would be:

\begin{center}
 \verb/AL.}AL.BNN.}D.OPGQO.}E.HRSJSUTUR.QGP.}!/
\end{center}
\medskip
\begin{center}
\includegraphics[scale=0.5]{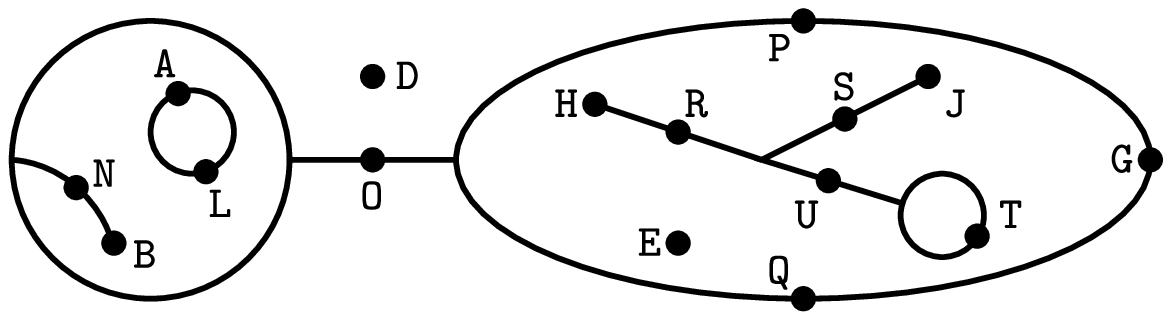}
\end{center}

\subsection{Generic vertices}

We replace the vertices that only occur once and in a boundary with a single vertex, by the generic vertex ``\verb/0/''. We also replace the vertices that only occur once but in a boundary with several vertices, by the generic vertex ``\verb/1/''.

The generic vertex ``\verb/2/'' can designate the vertices that occur twice in a row along a single boundary, just as \verb/N/ or \verb/O/ in our example. However, it can also designate the vertices that used to occur twice, and that now only occur once, because of a dead region deletion (such as \verb/T/ in our example).

\begin{center}
 \verb/AL.}AL.12.}0.2PGQ.}0.1RS1SU2UR.QGP.}!/
\end{center}
\medskip
\begin{center}
\includegraphics[scale=0.5]{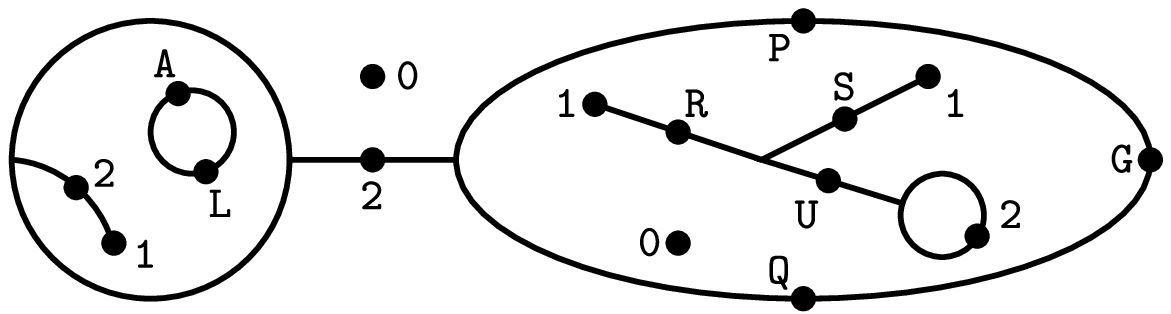}
\end{center}

\subsection{Lands}

We cut the position into independent parts, named \textit{lands}. A position can be cut into 2 lands when those lands have no letter in common.  The end-of-land character is ``\verb/]/''. Our example has 2 lands, and its string becomes:
\begin{center}
\verb/AL.}AL.12.}]0.2PGQ.}0.1RS1SU2UR.QGP.}]!/
\end{center}

\subsection{Renaming letters}

At this point, letters designate vertices that occur twice in the string, but we can distinguish between two different types: vertices that occur twice in the same boundary (which we will designate with lower-case letters) and vertices that appear in two different regions (which we will designate with upper-case letters).

We rename these vertices from ``\verb/a/'' and ``\verb/A/'', in the order of their appearance. We start again from ``\verb/a/'' when we meet a new boundary and we start again from ``\verb/A/'' when we meet a new land.

\begin{center}
\verb/AB.}AB.12.}]0.2ABC.}0.1ab1bc2ca.CBA.}]!/
\end{center}
\medskip
\begin{center}
\includegraphics[scale=0.5]{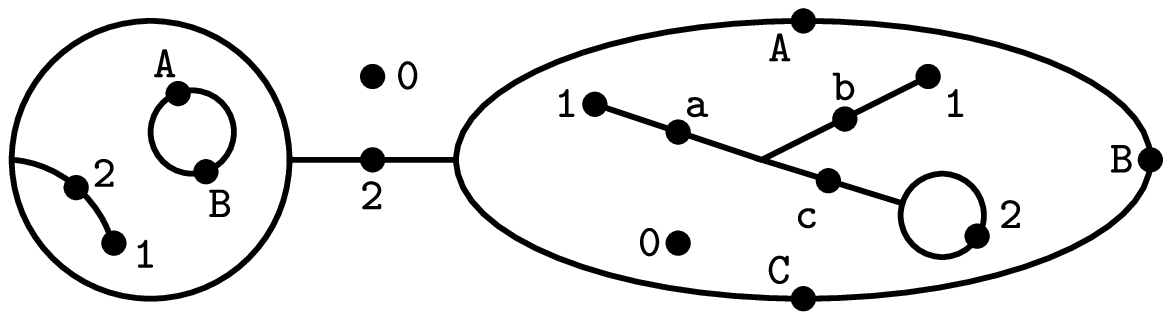}
\end{center}

\subsection{Regions equivalences}

When a region has 3 lives or less, we have an equivalent game if we merge the boundaries. For instance, the region \verb/A.BC.}/ becomes \verb/ABC.}/, or \verb/2.2.}/ becomes \verb/22.}/. These equivalences are a useful trick in our program, reducing considerably the number of stored strings.

\section{Strings canonization}
\label{canosec}

We already know that several strings can represent the same position, e.g. the previous position could also be represented by:
\begin{center}
\verb/BA2C.0.}0.2ca1ab1bc.CBA.}]AB.21.}AB.}]!/
\end{center}

As in \cite{ajs91}, we call \emph{canonization} the choice of a single string amongst all the equivalent strings. By merging equivalent branches in a game tree, canonization decreases efficiently both memory consumption and running time. On the other hand, canonization itself takes a lot of running time, and we need to take this into consideration when writing the program.

\subsection{Canonization}

First, we define an equivalence on the set of strings. Two strings are equivalent if they are equal modulo:

\begin{itemize}
\item the first vertex chosen in a boundary.
\item the boundaries order in a region.
\item the regions order in a land.
\item the lands order in the position.
\item the orientation chosen for each region.
\item a renaming of the vertices.
\end{itemize}

We can now define the term \emph{canonization}: this is the choice of a single string in each equivalent class. This choice should be as simple as possible, therefore, we will choose the minimal string for the following lexicographical order:
\begin{center}
\verb/0 < 1 < 2 < a < b < ... < A < B < ... < . < } < ] < !/
\end{center}

A little reflexion shows then that the canonized string of our previous position is \verb/0.1ab1bc2ca.ABC.}0.2ABC.}]12.AB.}AB.}]!/.

\begin{remark}
We are not allowed to change the orientation of a single boundary. We can only change the orientation of a complete region, i.e. the orientation of all the boundaries in the region.

For instance, \verb/122a2a.22.2AB.}2A.}2C.}BC.}]1122.}]!/ represents a losing position, while \verb/122a2a.22.2BA.}2A.}2C.}BC.}]1122.}]!/ represents a winning position\footnote{More obvious examples may exist, but this is the simplest that we know. However, to check this example, a long computation is necessary.}.
\end{remark}

\subsection{Pseudo-canonization}

We cannot actually perform the canonization, mainly because choosing the name of the upper-case letters requires too much running time. In particular, if a land has $k$ upper-case letters, performing the true canonization consists in choosing the minimal string for lexicographical order amongst the $k!$ possible strings.\par
Therefore, we only perform a pseudo-canonization: the same position can be represented by several strings, for instance the strings \verb/0.AB.CD.}0.AB.}CD.}]!/ and \verb/0.AB.CD.}0.CD.}AB.}]!/ represent the same position. Our pseudo-canonization algorithm renames the upper-case letters from ``\verb/A/'' in the order of their appearance and then sorts the string. Experience shows that performing this operation twice is the most efficient for both running time and memory usage. Ultimately, we only lose a few percentages of memory compared to a true canonization, whose time complexity prevents from computing $p$-spot games for $p\geq 5$ or $6$.\par
We have developed the complete game tree of the $p$-spot games, for $p\leq 6$, in order to evaluate the performance of our pseudo-canonization. Here is an extract of the complete game tree of the $4$-spot game:
\begin{center}
\begin{tabular}{ |c| }
\hline
\verb/ABCDEF.}ABCDEG.}FG.}]!/\tabularnewline
\verb/ABCDEF.}ABCDGF.}EG.}]!/\tabularnewline
\verb/ABCDEF.}ABCGEF.}DG.}]!/\tabularnewline
\verb/ABCDEF.}ABGDEF.}CG.}]!/\tabularnewline
\verb/ABCDEF.}BCDEFG.}AG.}]!/\tabularnewline
\hline
\end{tabular}
\end{center}

Since these strings represent the same position, a true canonization would have displayed only one string instead of these five. However, this position is not needed to compute the outcome of the $4$-spot game. In fact, the performance of the pseudo-canonization (comparatively to the true canonization) is better in a real computation that in a complete game tree development, because in a real computation we meet strings easier to compute, with less upper-case letters.\par
The following table gives the number of pseudo-canonized strings stored after a complete game tree development for the $p$-spot game.
\begin{center}
\begin{tabular}{ |c|c| }
\hline
n & number of strings\tabularnewline
\hline
2 & 18\tabularnewline
\hline
3 & 157\tabularnewline
\hline
4 & 1796\tabularnewline
\hline
5 & 24784\tabularnewline
\hline
6 & 393103\tabularnewline
\hline
\end{tabular}
\end{center}

This table could be useful for evaluating the performance of our canonization, in comparison with the canonization of other programs.

\section{Main algorithm}
\label{mainalgo}
\subsection{Outcome computation algorithm}

The standard recursive algorithm to compute the outcome of a position\footnote{In the following, we use the term ``position'' even if we speak of the string that represents it.} can be described as follows:

\begin{algorithm}
\label{cwl}
function compute-win-loss(position $P$)
\begin{itemize}
\item compute the children of $P$
\item for each child, do: if compute-win-loss(child)=Loss, return Win
\item return Loss\footnote{This step is reached only if no child is losing.}
\end{itemize}
\end{algorithm}

\begin{figure}[ht]
\centering
\includegraphics[width=.5\linewidth]{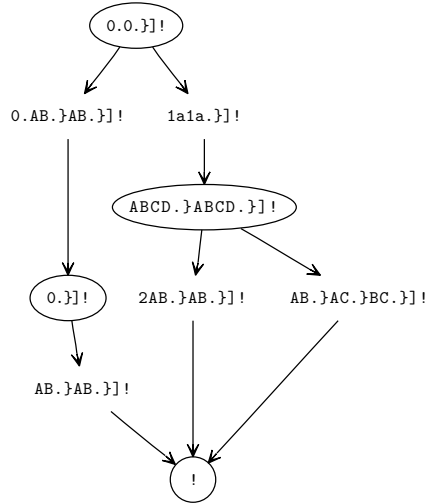}
\caption{2-spot game.}
\label{2spots}
\end{figure}

The core of our main algorithm uses this classical depth-first procedure. We also store the outcomes of the computed positions in a database (a \emph{transposition table}), in order not to compute several times the same position. Figure \ref{2spots} shows how this algorithm works with the $2$-spot game.

Surrounded positions are the losing ones. As figure \ref{2spots} shows, a position is losing if all its children are winning. But on the contrary, only one losing child is sufficient to prove that a position is winning, which implies that only a part of the complete game tree is sufficient to determine the outcome of the game (only 9 of the 18 positions of the complete game tree are needed here). We will detail this point again in section \ref{childrenOrdering}.

\subsection{Sum of independent games}

Sprouts positions can frequently be separated into independent games, the \emph{lands}. It would then be convenient to compute these lands independently, and deduce the outcome of the complete position from the outcome of the lands.

As explained in \cite{ajs91}, this can partially be done, because the sum of two losing games is losing, and the sum of a losing and a winning game is winning. However, the weakness of this method is that the sum of two winning games can either be winning or losing. We had implemented it in our program at first, and rather quickly, a large number of positions with several lands was saturating our databases, preventing us from improving the results of \cite{ajs91} of more than 2 or 3 spots.

As described below, the concept of \emph{nimber} can solve this problem, allowing us to compute the lands separately.

\subsection{Nimber theory}

A more detailed view of this theory can be found, amongst others, in \cite{ww01}.

\begin{definition}
The nimber of a position $P$ is denoted by $|P|$, and is defined as the smallest non-negative integer that is not a nimber of a child of $P$.
\end{definition}

This definition implies that $|P|=0$ if $P$ is losing, and $|P|\geq 1$ if $P$ is winning. We can also see that, if we know the complete game tree of $P$, we can recursively compute $|P|$, but in fact it is not necessary to expand the complete game tree to compute the nimber of a position, as the algorithm \ref{012} shows below.

The main result of the nimber theory can be stated as follows:

\begin{theorem}
If a position $P$ is made of two independent positions $P_1$ and $P_2$, then $|P|$ is the ``bitwise exclusive or'' of $|P_1|$ and $|P_2|$, denoted by $|P_1| \wedge |P_2|$.
\end{theorem}

For example, $|$\verb/1AB.}AB.}]!/$|=3$ and $|$\verb/22.}]!/$|=1$, so\\ $|$\verb/1AB.}AB.}]22.}]!/$|=3\wedge 1=2$.

\subsection{Couples}

In our program, instead of computing the outcome of a position, we compute the outcome of a \emph{couple}: (position+nimber), which represents the sum of two independent games. The position part consists of the game of Sprouts for the given position. The nimber part is the game of Nim for the given nimber value. The outcome of ($\emptyset$+$n$) is a loss if $n=0$, a win if $n\geq 1$. The original position $P_0$ is replaced by the couple ($P_0$+$0$).

We see that ``($P$+$n$) is losing'' means that $|P|=n$. As we only store losing positions in our program, it means that we will only store positions whose nimber is known (and the winning positions that we do not store correspond to positions whose nimber is only known to be different of certain values).

To determine the outcome of a couple, we can still use the algorithm \ref{cwl}, by extending to couples the definition of children: the children of a couple ($P$+$n$) are the children of the position part, whose form is (child($P$)+$n$), and the children of the nimber part, whose form is ($P$+$m$), with $m<n$.

\subsection{Computation of the nimber of a position}

With algorithm \ref{cwl} applied to a couple ($P$+$n$), we already are able to determine if $n$ is the nimber of the position $P$. If we need to compute the nimber $|P|$ of this position, we use the following (simple but efficient) method:

\begin{algorithm}
\label{012}
function nimber-of(position $P$)
\begin{itemize} 
 \item $n:=0$, $found:=false$
 \item while $found=false$ do:\par
if compute-win-loss($P$+$n$)=Loss, then $found:=true$, else $n:=n+1$
 \item return $n$
\end{itemize}

\end{algorithm}

It merely consists in trying $0$, $1$, $2$... until we find the right nimber. This algorithm will only be used on single lands, as explained below.

\subsection{Positions with several lands}

If a position is made of two lands, as in $(P_1]P_2]!+n)$, and if we know the nimber $|P_2|$, the above theorem shows that the outcome of $(P_1]P_2]!+n)$ will be the same as the outcome of $(P_1]!+n \wedge |P_2|)$.

So when our main algorithm meets a position with several lands, it computes with the algorithm \ref{012} the nimbers of all lands except one and merges the results with the nimber part of the couple. Therefore, we always compute the lands separately, and an indirect consequence is that we store only single lands in our database.

\subsection{Main algorithm}

With the previous ideas, algorithm \ref{cwl} is now:

\begin{algorithm}
\label{cwl2}
function compute-win-loss(couple($P$+$n$))
\begin{itemize}
\item for the lands of $P$ whose nimber is already stored in the database, merge their nimber with the nimber part (with bitwise xor).
\item compute the nimber of all the unknown lands of $P$ with algorithm \ref{012} (except for one land), and merge those nimbers with the nimber part ($P$ is then a single land).
\item for each child\footnote{Children of position part, and children of nimber part.} of ($P$+$n$), do:\par
	if compute-win-loss(child)=Loss, return Win
\item store ($P$+$n$) in the database and return Loss
\end{itemize}
\end{algorithm}

\begin{figure}[ht]
\centering
\includegraphics[width=.8\linewidth]{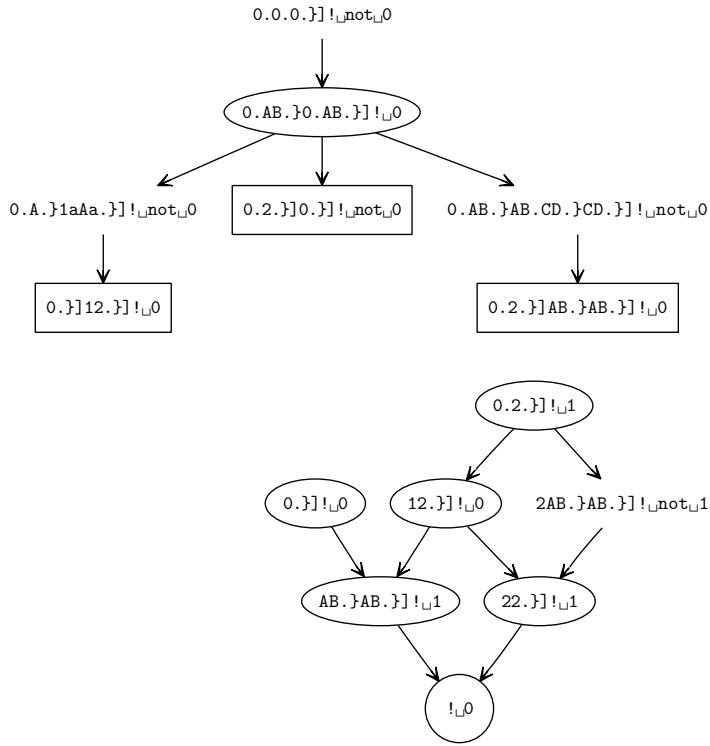}
\caption{3-spot game.}
\label{3spots}
\end{figure}

Figure \ref{3spots} shows how this algorithm works with the $3$-spot game.

Couples surrounded by a rectangle are those made of several lands. Their nimber is computed from the nimber of their lands, which can be found lower in the graph.

Couples surrounded by an ellipse are losing. Only those couples are stored by the program. All their children are computed to find their position's nimber.

The remaining couples are winning. The nimber of their position part is not known. We only know the nimber to be different of one (or sometimes several) value, because this value is already the nimber of some child of the position.

We see on this example the efficiency of this algorithm. There are 157 positions in the complete game tree of the $3$-spot game, but the program only needs 14 nodes to compute the outcome.

\section{Computation improvement}

\subsection{Positions storage}
\label{storing}
As in \cite{ajs91}, we store only losing couples in order to reduce the memory consumption: as the complexity of the computation increases, losing couples become less frequent, so this choice allows us to reduce the size of the transposition table. The profit is variable, but we can take into consideration that the size of the database is divided by a factor between 5 and 20, which is not negligible. In exchange, we need a little more running time when computing again a winning couple that has already been computed before. Since the couple is not stored in the database, its winning outcome will be known only after a computation of its children, one of which will be found to be losing in the database. \par
The main difference of our database compared to \cite{ajs91} is that we store losing couples (positions whose nimber is known), and not losing positions. Moreover, as explained before, we store only positions with a single land.\par
It is also possible in our program to export databases as text files, allowing us to analyze them manually, or resume later an incomplete computation.

\subsection{Children ordering}
\label{childrenOrdering}
When a position (or a couple) is losing, there is no short way to prove it: we need to compute all the children and prove that they are all winning. On the contrary, when a position is winning, we only need to find one losing child, so that there are several ways to search the game tree and compute the outcome of the starting position. These ways are not of equivalent difficulty, because Sprouts game trees are unbalanced, with some areas far more complicated than others.

Consequently we sort the children by computational difficulty, in order to find a losing child as quickly as possible. Thus, we lose the least possible time possible in unnecessary computations of winning children and also find first the losing children whose outcomes are easy to compute.

When ordering the children, we try to evaluate their difficulty only from their string. The rules used to evaluate the difficulty and order the children are an important part of the program, because a bad set of rules could push the computation into complicated parts of the game tree. Our current set of rules is as follows:

\begin{itemize}
\item priority to couples with minimal (number of lives + nimber).
\item priority to positions with a lot of lands.
\item priority to positions with a little estimated number of children.
\end{itemize}

Since it would require too much running time to compute the exact value, the number of children is only estimated as follows:

\begin{itemize}
\item When linking a boundary to itself, the possible number of children is $(number~of~vertices)^2\times (number~of~partitions)$, where
\\ $(number~of~partitions)$ is the number of ways to partition the other boundaries into two sets.
\item When linking two different boundaries, the possible number of children is the sum, for any couple $(B_i, B_j)$ of boundaries, of\\ $(number~of~vertices~of~B_i)\times(number~of~vertices~of~B_j)$.\par
\end{itemize}

We must also make a choice when considering a position composed of several lands. We decided to compute first the nimbers of the lands with the smallest number of lives.

These rules significantly improve the computation compared to a random exploration of the game tree, but they are still far from providing an optimal exploration of the game tree. In fact, we believe that whatever the rules used to order the children, their efficiency is eventually limited and more global search algorithms are needed.

\subsection{Manual exploration of the game tree}

We have implemented an interface to follow and interact in real-time with the computation process. When the computation seems to be stuck in some part of the game tree, we can manually decide to explore elsewhere.

During the computation, the program displays the list of currently studied couples. At a given instant, the program is computing the outcome of a list of couples, each one being the child of the previous one. For example, here is what the program could display during a $12$-spot game computation:
\begin{center}
\begin{tabular}{ |c|c|c| }
\hline
level & position part & nimber part\tabularnewline
\hline
1 & \verb/0.0.0.0.0.0.AB.}0.0.0.0.0.AB.}]!/ & 0\tabularnewline
2 & \verb/0.0.0.0.0.}]0.0.0.A.}0.0.0.A.}]!/ & 0\tabularnewline
3 & \verb/0.0.0.0.0.}]!/ & 1\tabularnewline
4 & \verb/0.0.0.0.AB.}AB.}]!/ & 1\tabularnewline
5 & \verb/0.0.0.0.}]!/ & 1\tabularnewline
6 & \verb/0.0.1a1a.}]!/ & 1\tabularnewline
\hline
\end{tabular}
\end{center}

With a click in the interface, it is then possible to choose another child on a given level, and so to compute another part of the game tree. For example, we could decide to compute the couple (\verb/0.0.A.}0.0.A.}]!/+1) in the fifth level.

When the position is composed of several lands, the program computes the nimber of one of these lands with algorithm \ref{012} (here, it has already computed that the nimber of \verb/0.0.0.0.0.}]!/ is not 0, and is now trying with 1). In this case too, we can click in the interface to compute first the nimber of another land. In our example, we could decide to compute the nimber of \verb/0.0.0.A.}0.0.0.A.}]!/ in the third level.

We empirically decide where and when we should click: if the program spends too much time on a branch of the game tree, we decide to change it. If a couple on a level is almost computed to be losing (which means that almost all its children have been computed to be winning), it is often efficient to click two levels below, to try to quickly find losing children of the remaining unknown children. We can also have a look at the position to try to choose easier couples (positions are easier when they tend to be quickly cut into lands).

Being able to manually choose which part of the game tree to explore was far more efficient that any automatic selection that we imagined. For example, whereas a computation for the $12$-spot game ends after having stored more than $100,000$ couples without human intervention, by clicking, somebody with some experience can end this game in less than $2,000$ couples.

\subsection{Check computation}

User interactions proved to be powerful, but they have a drawback: it is impossible to reproduce exactly the same computation twice. For this reason, when a computation succeeds, we perform a standalone check computation, which uses the previous results to guide itself in the game tree. This check computation does not authorize any interaction from the user and is of course reproducible for a given database of previous results.

The check computation also reduces the number of couples needed to prove the result. Indeed, during the first computation, we compute many useless couples (the winning children of winning couples are usually useless). When  we meet a new couple in the check computation, we look for its value in the database of previous results, and if it is losing (i.e. a position whose nimber is known), we compute again all its children to check it. But if it is winning, we look in the previous database which child is losing, and we compute only this one.

The result of a check computation is then a \textit{solution tree}, providing a winning strategy without any unnecessary information. These solution trees are small enough to provide hand-checkable proofs for $p$-spot games with little values of $p$. For bigger values of $p$, it is mainly a way of reducing the size of the files storing the winning strategies.

\begin{figure}[ht]
\centering
\includegraphics[width=\linewidth]{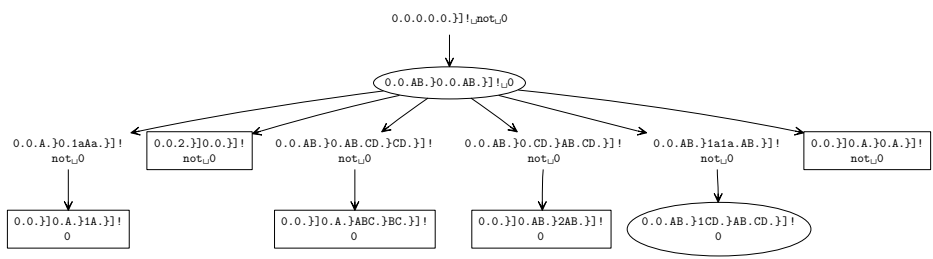}
\caption{A solution for the 5-spot game, positions with 12 lives or more.}
\label{5spots}
\end{figure}

During the check computation, we can generate a file compatible with the graph visualization tool Graphviz\footnote{http://www.graphviz.org/}, which enables us to draw graphs describing how the computation proceeds (e.g. figures \ref{2spots} and \ref{3spots}). For bigger values of $p$, the complete graph would be rather unreadable, but we can choose to plot only the positions with a minimal number of lives. For example, figure \ref{5spots} shows the upper part of the graph for the 5-spot game, with only the positions with 12 lives or more.

Here, we can see how cutting positions into lands simplifies the computation: to complete it, we only need to compute the nimber of 6 lands of less than 7 lives, and the outcome of a position of 12 lives. The most spectacular example is the case of $p$=17: in our computation, every position is cut into several lands of less than 27 lives after only 3 moves from the starting position (which has 51 lives).

\begin{figure}[ht]
\centering
\includegraphics[width=\linewidth]{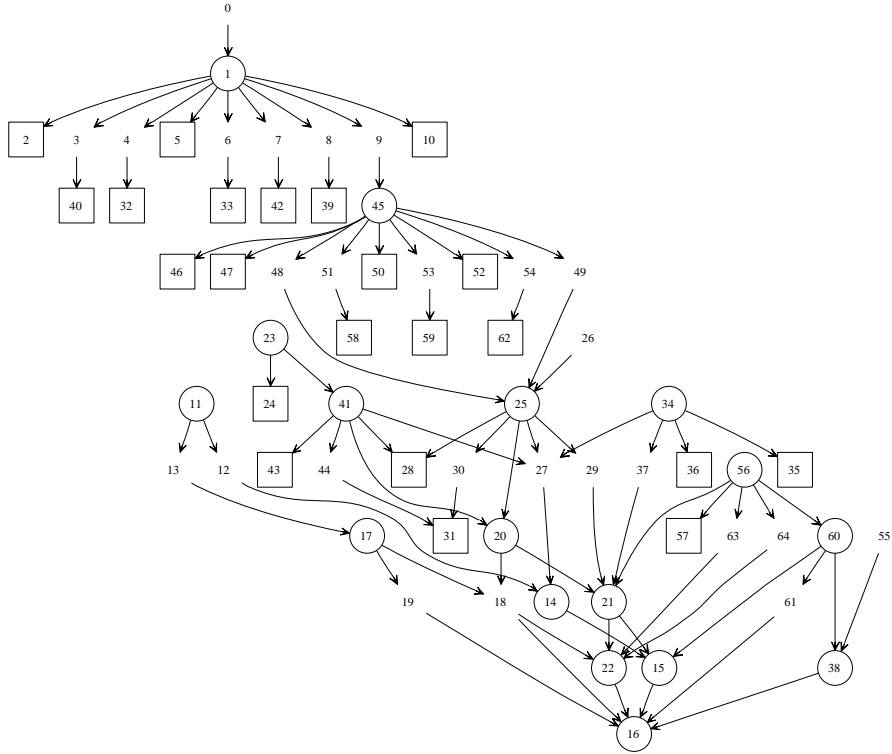}
\caption{A complete solution for the 4-spot game.}
\label{4spots}
\end{figure}

It is also possible to plot only reference numbers instead of the whole positions, which is useful for bigger graphs. See, for instance, figure \ref{4spots}, which shows the complete 4-spot game. The numbers have no special meaning, the positions which they refer to are listed in the appendix, with their computed values.

\section{Results}

\subsection{$p$-spot games already computed}
\label{nspotres}
The following table shows the number of couples stored in the minimal databases for all the $p$-spot games computed up to now. We have computed all $p$-spot games up to 32 spots. The 33-spot game is the first unknown one, but the biggest computed one is the 47-spot game. The results so far all support the conjecture emitted in \cite{ajs91}: the first player loses if and only if $p$ is $0$, $1$ or $2$ modulo $6$.

\begin{center}
\begin{minipage}{0.13\textwidth}
\begin{tabular}{ |c|c| }
\hline
$p$ & size\tabularnewline
\hline
0 & 0\tabularnewline
1 & 1\tabularnewline
2 & 3\tabularnewline
3 & 6\tabularnewline
4 & 16\tabularnewline
5 & 38\tabularnewline
6 & 64\tabularnewline
\hline
\end{tabular}
\end{minipage}
\begin{minipage}{0.14\textwidth}
\begin{tabular}{ |c|c| }
\hline
$p$ & size\tabularnewline
\hline
7 & 103\tabularnewline
8 & 205\tabularnewline
9 & 63\tabularnewline
10 & 140\tabularnewline
11 & 140\tabularnewline
12 & 475\tabularnewline
13 & 577\tabularnewline
\hline
\end{tabular}
\end{minipage}
\begin{minipage}{0.15\textwidth}
\begin{tabular}{ |c|c| }
\hline
$p$ & size\tabularnewline
\hline
14 & 1580\tabularnewline
15 & 3252\tabularnewline
16 & 1068\tabularnewline
17 & 471\tabularnewline
18 & 3233\tabularnewline
19 & 3630\tabularnewline
20 & 4051\tabularnewline
\hline
\end{tabular}
\end{minipage}
\begin{minipage}{0.16\textwidth}
\begin{tabular}{ |c|c| }
\hline
$p$ & size\tabularnewline
\hline
21 & 9270\tabularnewline
22 & 5706\tabularnewline
23 & 2837\tabularnewline
24 & 9316\tabularnewline
25 & 9229\tabularnewline
26 & 18567\tabularnewline
27 & 59117\tabularnewline
\hline
\end{tabular}
\end{minipage}
\begin{minipage}{0.16\textwidth}
\begin{tabular}{ |c|c| }
\hline
$p$ & size\tabularnewline
\hline
28 & 14813\tabularnewline
29 & 3414\tabularnewline
30 & 58363\tabularnewline
31 & 58365\tabularnewline
32 & 58204\tabularnewline
33 & ?\tabularnewline
34 & 21107\tabularnewline
\hline
\end{tabular}
\end{minipage}
\begin{minipage}{0.16\textwidth}
\begin{tabular}{ |c|c| }
\hline
$p$ & size\tabularnewline
\hline
35 & 18812\tabularnewline
... & ?\tabularnewline
40 & 45782\tabularnewline
41 & 48890\tabularnewline
... & ?\tabularnewline
47 & 54542\tabularnewline
... & ?\tabularnewline
\hline
\end{tabular}
\end{minipage}
\end{center}

It is rather surprising that the number of stored couples does not increase strictly with $p$. In fact, some patterns seem to occur modulo $6$.

We found that $p=15$, $21$ and $27$ were very difficult to compute, because they are winning, and only have one losing child, which is the most complicated one (its string is \verb/0.0.0.(...)0.1a1a.}]!/). $p=33$, the first value that we do not yet know, seems to follow the same path.

Conversely, $p=17$, $23$, $29$, $35$, $41$ and $47$ were much easier to compute. They are winning positions, and the losing child is obtained by linking one spot to itself and separating the remaining spots in two equal sets (this child is often the easiest to compute).

We also noticed that the number of stored couples is almost equal for $p=15$ and $p=18,19$: once the result of $p=15$ is known, we only needed a little more computation to deduce the result for $p=18,19$. The same phenomenon occurs for $p=21$ and $p=24,25$, or for $p=27$ and $p=30,31$.

We tried to minimize the number of couples for $p\leq 25$ (and $p=29$), but we did not take the time to do the same work for the other values, therefore the reader needs to be aware that it is probably possible to reduce significantly the number of couples for those values.

\subsection{Nimber conjectures}

We observed that the nimber of the winning starting positions, for $p\leq 32$, is $1$, so we propose a stronger conjecture than the ``Sprouts conjecture'' emitted in \cite{ajs91}:

\begin{conjecture}
The nimber for the starting position with $p$ spots is $0$ if $p$ is $0$, $1$ or $2$ modulo $6$, and $1$ if $p$ is $3$, $4$ or $5$ modulo $6$.
\end{conjecture}

It is rather easy to imagine other conjectures around the nimber. For example, if we compute the nimber of the position: \verb/222(...)222.}]!/, with $p$ generic vertices ``\verb/2/'', we obtain (starting from \verb/22.}]!/): 1;0;2;1;0;3;1;0;5;1;0;3;1;0;3;1;0;3;1;0;3;1;0.

We also imagined another extension of the Sprouts conjecture: the nimber of the position is left unchanged by an addition of 6 boundaries ``\verb/0./'' into a given region. This conjecture is false, since \verb/0.22.}]!/ has nimber 0, and \verb/0.0.0.0.0.0.0.22.}]!/ has nimber 2. But it does work for more than 90\% of the positions that we tried. The existence of such nearly-true patterns tends to infirm the Sprouts conjecture.

\subsection{Hand-checkable proofs}

Using a computer to determine mathematical results is not completely convincing, especially because there could be a programming error, considering the size of the program. So, the most skeptic readers could use the program to generate files that would provide them with hand-checkable proofs of the simplest results. For example, figure \ref{4spots} provides a proof that the first player wins in the 4-spot game.

We printed the result of a check computation for the 9-spot game and manually checked 66 losing couples (positions of known nimber), and 258 winning couples (positions whose nimber is known to be different of certain values). For the losing couples, we checked that we obtained the same sets of children as our program, and for the winning ones, we had to check only one child, very seldom two or more. This took us a few hours.

The table in paragraph \ref{nspotres} implies that we could use this method to check in a reasonable amount of time the value of the $p$-spot game with $p\leq 11$. The program can also point in the right direction someone who would like to create a totally manual proof, like the one in \cite{fl04}.

\section*{Conclusion}

The main obstacle for computing higher $p$-spot games is neither memory nor computation time, but rather human time for the manual exploration of the game tree. This human intervention is an embryo of a best-first search, which is probably more suited for the game of Sprouts than the depth-first search currently performed by our program. Therefore, implementing classical best-first search algorithms, such as the PN-search, would probably lead to better results.
\par
Distributed computing is also another solution to compute higher $p$-spot games, and the check computation, by reducing the size of the databases, would be really efficient for this. It is likely that in the coming years, we will know who wins in the game of Sprouts when starting with more than fifty spots.

\medskip

The program that we used for the computations is available with its source code on our web site \verb|http://sprouts.tuxfamily.org/| under a GNU licence, together with several resulting databases.

\section*{Acknowledgements}

We wish to thank Jean-Paul Delahaye for his decisive contribution to the continuation of our work.

%------------------------------------------------------------------------------
% Bibliography
%------------------------------------------------------------------------------

\bibliography{sprouts}
\bibliographystyle{amsplain}

\newpage

\section*{Appendix: Correspondence for Figure \ref{4spots}}

\begin{table}[ht]
\begin{minipage}[t]{.505\linewidth}
\begin{tabular}{ |c|c|c| }
\hline
\# & position & nimber \tabularnewline
\hline
0 & \verb/0.0.0.0.}]!/ & $\neq$ 0\tabularnewline
\hline
1 & \verb/0.0.AB.}0.AB.}]!/ & 0\tabularnewline
\hline
2 & \verb/0.0.2.}]0.}]!/ & $\neq$ 0\tabularnewline
\hline
3 & \verb/0.0.AB.}AB.CD.}CD.}]!/ & $\neq$ 0\tabularnewline
\hline
4 & \verb/0.0.A.}1aAa.}]!/ & $\neq$ 0\tabularnewline
\hline
5 & \verb/0.0.}]0.2.}]!/ & $\neq$ 0\tabularnewline
\hline
6 & \verb/0.1aAa.}0.A.}]!/ & $\neq$ 0\tabularnewline
\hline
7 & \verb/0.AB.CD.}0.AB.}CD.}]!/ & $\neq$ 0\tabularnewline
\hline
8 & \verb/0.AB.}0.CD.}AB.CD.}]!/ & $\neq$ 0\tabularnewline
\hline
9 & \verb/0.AB.}1a1a.AB.}]!/ & $\neq$ 0\tabularnewline
\hline
10 & \verb/0.A.}0.A.}]0.}]!/ & $\neq$ 0\tabularnewline
\hline
11 & \verb/0.0.}]!/ & 0\tabularnewline
\hline
12 & \verb/0.AB.}AB.}]!/ & $\neq$ 0\tabularnewline
\hline
13 & \verb/1a1a.}]!/ & $\neq$ 0\tabularnewline
\hline
14 & \verb/0.}]!/ & 0\tabularnewline
\hline
15 & \verb/AB.}AB.}]!/ & 1\tabularnewline
\hline
16 & \verb/!/ & 0\tabularnewline
\hline
17 & \verb/ABCD.}ABCD.}]!/ & 0\tabularnewline
\hline
18 & \verb/2AB.}AB.}]!/ & $\neq$ 0 1\tabularnewline
\hline
19 & \verb/AB.}AC.}BC.}]!/ & $\neq$ 0\tabularnewline
\hline
20 & \verb/0.2.}]!/ & 1\tabularnewline
\hline
21 & \verb/12.}]!/ & 0\tabularnewline
\hline
22 & \verb/22.}]!/ & 1\tabularnewline
\hline
23 & \verb/0.A.}0.A.}]!/ & 1\tabularnewline
\hline
24 & \verb/0.}]12.}]!/ & 0\tabularnewline
\hline
25 & \verb/0.AB.}2AB.}]!/ & 0\tabularnewline
\hline
26 & \verb/0.0.2.}]!/ & $\neq$ 0\tabularnewline
\hline
27 & \verb/0.A.}2A.}]!/ & $\neq$ 0\tabularnewline
\hline
28 & \verb/0.}]22.}]!/ & $\neq$ 0\tabularnewline
\hline
29 & \verb/1aAa.}2A.}]!/ & $\neq$ 0\tabularnewline
\hline
30 & \verb/2AB.}AB.CD.}CD.}]!/ & $\neq$ 0\tabularnewline
\hline
31 & \verb/22.}]AB.}AB.}]!/ & 0\tabularnewline
\hline
32 & \verb/0.0.}]12.}]!/ & 0\tabularnewline
\hline
\end{tabular}
\end{minipage}
\hfill
\begin{minipage}[t]{.485\linewidth}
\begin{tabular}{ |c|c|c| }
\hline
\# & position & nim. \tabularnewline
\hline
33 & \verb/0.A.}1A.}]0.}]!/ & 0\tabularnewline
\hline
34 & \verb/0.A.}1A.}]!/ & 0\tabularnewline
\hline
35 & \verb/0.}]AB.}AB.}]!/ & $\neq$ 0\tabularnewline
\hline
36 & \verb/12.}]1.}]!/ & $\neq$ 0\tabularnewline
\hline
37 & \verb/1A.}ABC.}BC.}]!/ & $\neq$ 0\tabularnewline
\hline
38 & \verb/1.}]!/ & 1\tabularnewline
\hline
39 & \verb/0.AB.}2AB.}]0.}]!/ & 0\tabularnewline
\hline
40 & \verb/0.A.}0.A.}]AB.}AB.}]!/ & 0\tabularnewline
\hline
41 & \verb/0.A.}ABC.}BC.}]!/ & 0\tabularnewline
\hline
42 & \verb/0.A.}ABC.}BC.}]0.}]!/ & 0\tabularnewline
\hline
43 & \verb/12.}]AB.}AB.}]!/ & $\neq$ 0\tabularnewline
\hline
44 & \verb/ABC.}ADE.}BC.}DE.}]!/ & $\neq$ 0\tabularnewline
\hline
45 & \verb/0.AB.}1CD.}AB.CD.}]!/ & 0\tabularnewline
\hline
46 & \verb/0.2.}]1AB.}AB.}]!/ & $\neq$ 0\tabularnewline
\hline
47 & \verb/0.AB.}2AB.}]1.}]!/ & $\neq$ 0\tabularnewline
\hline
48 & \verb/0.AB.}2CD.}AB.CD.}]!/ & $\neq$ 0\tabularnewline
\hline
49 & \verb/0.AB.}ABC.}CDE.}DE.}]!/ & $\neq$ 0\tabularnewline
\hline
50 & \verb/0.AB.}AB.}]12.}]!/ & $\neq$ 0\tabularnewline
\hline
51 & \verb/0.A.}1B.}aAaB.}]!/ & $\neq$ 0\tabularnewline
\hline
52 & \verb/0.}]1AB.}2AB.}]!/ & $\neq$ 0\tabularnewline
\hline
53 & \verb/1aAa.}1BC.}ABC.}]!/ & $\neq$ 0\tabularnewline
\hline
54 & \verb/1AB.}AB.CD.}/ & $\neq$ 0\tabularnewline
& \verb/CD.EF.}EF.}]!/ & \tabularnewline
\hline
55 & \verb/1AB.}AB.}]!/ & $\neq$ 1\tabularnewline
\hline
56 & \verb/1AB.}2AB.}]!/ & 1\tabularnewline
\hline
57 & \verb/1.}]22.}]!/ & 0\tabularnewline
\hline
58 & \verb/0.}]1.}]AB.}AB.}]!/ & 0\tabularnewline
\hline
59 & \verb/12.}]1A.}2A.}]!/ & 0\tabularnewline
\hline
60 & \verb/1A.}2A.}]!/ & 0\tabularnewline
\hline
61 & \verb/2A.}2A.}]!/ & $\neq$ 0\tabularnewline
\hline
62 & \verb/1AB.}2AB.}]AB.}AB.}]!/ & 0\tabularnewline
\hline
63 & \verb/2AB.}2AB.}]!/ & $\neq$ 1\tabularnewline
\hline
64 & \verb/2A.}ABC.}BC.}]!/ & $\neq$ 1\tabularnewline
\hline
\end{tabular}
\end{minipage}
\end{table}

\end{document}